\documentclass[12pt]{article}
\usepackage{fullpage,graphicx,psfrag,amsmath,verbatim,xcolor}
\usepackage[small,bf]{caption}

\usepackage{amssymb}
\usepackage{enumitem}
\usepackage{listings}
\usepackage{authblk}
\usepackage{textcomp}

\usepackage[noend]{algpseudocode}
\usepackage{algorithm}
\usepackage{url}
\usepackage{bbm}

\usepackage{hyperref}
\hypersetup{ colorlinks   = true,    
urlcolor     = blue,    
linkcolor    = blue,    
citecolor    = red      
}

\usepackage{tikz}
\usetikzlibrary{shapes.geometric}
\usetikzlibrary{positioning}

\usepackage{caption}
\usepackage{subcaption}
\usepackage{multirow}

\newcommand{\ones}{\mathbf 1}
\newcommand{\reals}{{\mbox{\bf R}}}

\newcommand{\symm}{{\mbox{\bf S}}}  

\newcommand{\BEQ}{\begin{equation}}
\newcommand{\EEQ}{\end{equation}}
\newcommand{\BEAS}{\begin{eqnarray*}}
\newcommand{\EEAS}{\end{eqnarray*}}

\newcommand{\BBM}{\left[\begin{matrix}}
\newcommand{\EBM}{\end{matrix}\right]}

\newcommand{\BIT}{\begin{itemize}}
\newcommand{\EIT}{\end{itemize}}
\newcommand{\BNUM}{\begin{enumerate}}
\newcommand{\ENUM}{\end{enumerate}}

\newcommand{\Tr}{\mathop{\bf Tr}}
\newcommand{\diag}{\mathop{\bf diag}}

\newcommand{\Expect}{\mathop{\bf E{}}}

\newcommand{\argmin}{\mathop{\rm argmin}}

\newcommand{\eg}{{\it e.g.}}
\newcommand{\ie}{{\it i.e.}}

\newcommand{\PAR}[1]{\paragraph{#1.}}

\title{
Value-Gradient Iteration\\ with Quadratic Approximate Value Functions
}
\author{Alan Yang \quad Stephen Boyd}

\graphicspath{ {./figures/} }

\begin{document}
\maketitle

\begin{abstract}
We propose a method for designing policies for convex stochastic control
problems characterized by random linear dynamics and convex stage cost. We
consider policies that employ quadratic approximate value functions as a
substitute for the true value function. Evaluating the associated control policy
involves solving a convex problem, typically a quadratic program, which can be
carried out reliably in real-time. Such policies often perform well even when
the approximate value function is not a particularly good approximation of the
true value function. We propose value-gradient iteration, which fits the
gradient of value function, with regularization that can include constraints
reflecting known bounds on the true value function. Our value-gradient iteration
method can yield a good approximate value function with few samples, and little
hyperparameter tuning. We find that the method can find a good policy with
computational effort comparable to that required to just evaluate a control
policy via simulation.

\end{abstract}

\newpage
\tableofcontents
\newpage

\section{Introduction}

We consider convex approximate dynamic programming (ADP) policies for convex
stochastic control problems, which involve systems with known random linear
dynamics and convex stage costs. Evaluating an ADP policy reduces to solving a
convex optimization problem involving a convex approximate value function. We
focus on fitting quadratic approximate value functions, and refer to the
associated policies as quadratic approximate dynamic programming (QADP) policies.
While QADP policies are optimal for problems with convex quadratic stage cost
\cite{Bertsekas2012,BarrattB2021}, they can also serve as effective heuristics
for other problem types. It has been observed that ADP policies can perform well
even when using imperfect approximations of the true value function
\cite{Powell2007,KeshavarzB2014,Bertsekas2019}.

In this work, we propose an approximate value iteration method for finding
quadratic approximate value functions for convex stochastic control problems,
which we refer to as \emph{value-gradient iteration} (VGI). 
In principle, an optimal
value function may be found by iterating the Bellman operator, which maps
real-valued functions on the state space to real-valued functions on the state
space \cite{Bellman1954}. Since it is not possible in general to exactly
represent functions on $\reals^n$, we incorporate a function approximation step
after each application of the Bellman operator, a general approach called
fitted value iteration (FVI).  In our proposed VGI, instead of directly
fitting the value function, we fit the gradient of the value function with
respect to the state vector.

It is sufficient to approximate the gradient since constant offsets in the
value function have no impact on the associated ADP policy. In addition, the
gradient of the value function carries more information than the value function
itself \cite{DayanS1995,Fairbank2008}. If the gradient is well approximated at a
set of states, then the value function is also well approximated locally around
those states, up to an additive constant which does not affect the policy.
However, having a good approximation
of only the value at a set of states does not imply that the value function is
well approximated locally around those states.

Most importantly, VGI is practical to implement for QADP. We show that, when it
exists, the gradient of the Bellman operator applied to a convex quadratic
function can be obtained at any state by evaluating a particular optimal dual
variable associated with the QADP policy. Since the gradient of a convex
quadratic is an affine function, in each iteration we fit an affine function to
a set of pairs of states and value-gradients. This fitting problem is a convex
optimization problem. Therefore, VGI involves solving a sequence of convex
optimization problems, which can be carried out reliably.

We also consider several techniques for enhancing the reliability of VGI,
including damping, a robust Huber fitting loss, and the incorporation of prior
knowledge constraints and regularization. VGI remains effective even when the
state space dimension is large relative to the number of fitting samples, as we
will demonstrate with several numerical examples. Finally, we note that the
computational effort of obtaining a good QADP policy using VGI is small enough
that it is comparable to that of simply evaluating the policy through
simulation.

\subsection{Related work}

\PAR{Dynamic programming}
Dynamic programming (DP) provides techniques for computing the optimal value
function and policy for general Markov decision processes. The optimal policy is
evaluated by solving an optimization problem, where the control is chosen by
minimizing the current stage cost plus the expected value function at the next
state. For convex stochastic control problems, this is a convex optimization
problem \cite{Bellman1954,BertsekasS1996,Bertsekas2017,Puterman2014}. However,
it is possible to exactly represent and find the value function in a only 
few special
cases, for example when the state space is discrete \cite{SuttonB2018}, or when
we have a convex stochastic control problem with a convex extended quadratic
stage cost \cite{BarrattB2021}.

\PAR{Approximate dynamic programming}
ADP \cite{Powell2007,Bertsekas2012,Bertsekas2019} methods are heuristics used in
stochastic control when the problem cannot be solved by applying DP directly.
Typically, these methods either approximate the value function in DP or tune the
parameters of a parametric policy. In some contexts, approximate value functions
are known as control Lyapunov functions \cite{FreemanP1996, CorlessL1998}.

One approach to ADP is to approximate the value function by relaxing the Bellman
equation to an inequality, and then solving a convex optimization problem
involving a model of the dynamics and stage cost. When the state and input
spaces are finite, this leads to a linear program (LP) \cite{DeFariasVR2003}.
When the dynamics are affine, the stage cost is quadratic, and the input is
constrained to be in a convex set, quadratic approximate value functions can be
obtained using semidefinite programming \cite{WangB2009,WangOB2015}. In both
cases, the resulting approximate value functions are lower bounds on the true
value function.

Other value function approximation methods search for an approximate value
function that satisfies the Bellman equation along simulated trajectories. This
includes the method proposed in this paper, which is closely related to fitted
(or projected) value iteration \cite{BellmanD1959, KeshavarzB2014,
Bertsekas2012}. Other methods, which do not assume that a model of the dynamics
and stage cost are available, include $Q$-iteration
\cite{AntosSM2007,Bertsekas2012}, $Q$-learning \cite{WatkinsD1992,SuttonB2018},
and temporal difference learning \cite{Sutton1998, BertsekasBN2004}.

Instead of approximating the value function, other ADP techniques directly
optimize the parameters of a parametric policy to improve performance along
system trajectories. Stochastic gradient descent and its variants have been
employed to tune convex optimization control policies \cite{AgrawalBBS2020} and
controllers based on Proportional-Integral-Derivative (PID) control
\cite{Minorsky1922,AastromHHH1993} and model predictive control (MPC)
\cite{CamachoA2013,AmosJSBK2018}. Policy gradient methods provide a method for
differentiating through policies parametrized by neural networks
\cite{MnihBMGLHSK2016,SchulmanWDRK2017}. 

\PAR{Reinforcement learning}
Reinforcement learning (RL) methods \cite{SuttonB2018,Bertsekas2019} can be
considered a form of approximate dynamic programming (ADP), although their
primary focus is on learning from interactions with the system or a simulator,
rather than relying on explicit mathematical models of the system dynamics or
stage cost. In this work, we assume that models of the dynamics and stage cost
are either known or have been estimated or learned beforehand. This is similar
to some model-based RL methods that learn a policy and a model of the dynamics
jointly \cite{Sutton1990,DeisenrothR2011}. In the context of control, the
process of learning the dynamics is typically referred to as system
identification \cite{Ljung1998}. 

\PAR{Value gradients}
When considering a differentiable approximate value function, it is advantageous
to have accurate approximations of its derivatives with respect to the state,
\ie, the value gradient. If the value gradient is well-approximated along a
simulated trajectory, then the approximate value function also provides a good
local approximation around that trajectory \cite{DayanS1995}. Notably, it is
only necessary to approximate the value gradient since constant offsets in the
approximate value function do not affect the associated policy.

On the other hand, solely having a good approximation of the value function
itself along a trajectory does not ensure a good local approximation. In many
cases, value function approximation methods rely on stochastic local
exploration, such as dithering \cite{Bertsekas2012, SuttonB2018}, to overcome
this limitation. Indeed, value-gradient-based RL methods such as dual heuristic
programming (DHP), \cite{Werbos1999} globalized DHP \cite{ProkhorovW1997},
value-gradient learning \cite{Fairbank2008,FairbankA2012}, and stochastic value
gradients \cite{HeessWSLET015,AmosSYW2021} have been shown to find better
policies using less simulation than value function approximation methods that do
not directly approximate the value gradient.

VGI differs from the aforementioned value-gradient-based methods in that it does
not require stochastic approximations of the value gradient. Fitted value
iteration with value gradients is tractable for convex stochastic control
problems, since we can exactly evaluate the gradient of the Bellman operator
applied to a convex approximate value function by solving a convex optimization
problem.

\PAR{Convex optimization control policies}
For convex stochastic control, the policy associated with a convex quadratic
approximate value function can be evaluated by solving a convex optimization
problem, \ie, it is a convex optimization control policy (COCP)
\cite{AgrawalBBS2020}. COCPs are typically evaluated by solving quadratic
programs (QPs), which can often be done efficiently in real-time
\cite{WangB2010}. Evaluating a COCP may also involve minimizing a more complex
convex function, such as one parametrized by a neural network \cite{AmosXK2017}
To enable embedded applications, code generation tools like CVXGEN
\cite{MattingleyB2012} and CVXPYgen \cite{SchallerBDASB2022} can be utilized.

Other examples of COCPs include convex model predictive control (MPC)
\cite{GarciaPM1989,BorrelliBM2017} and convex approximate dynamic programming
\cite{KeshavarzB2014,Bertsekas2012}. COCPs can also be tuned by differentiating
through their solution maps \cite{AgrawalBBS2020,AmosJSBK2018}.

\subsection{Outline}

In \S\ref{s-cvx-stoch-control}, we introduce the convex stochastic control
problem and solution methods, via dynamic programming and model predictive
control. Approximate dynamic programming with quadratic approximate value
functions is described in \S\ref{s-qadp}, value-gradient iteration is introduced
in \S\ref{s-qfvgi}, and extensions and variations are discussed in
\S\ref{s-extensions}. In \S\ref{s-examples}, we present three numerical
examples: an input-constrained linear quadratic regulator (LQR) problem, a
commitments planning problem involving alternative investments, and a supply
chain optimization problem.

\section{Convex stochastic control}\label{s-cvx-stoch-control}

\subsection{Average-cost convex stochastic control problem}

\PAR{Dynamics} We consider a dynamical system evolving in discrete time
$t=0,1,2, \ldots$, with state $x_t\in\reals^n$, input $u_t\in\reals^m$, and
affine dynamics
\[
x_{t+1} = A_tx_t + B_t u_t + c_t, \quad t=0,1, \ldots,
\]
where $A_t\in\reals^{n\times n}$, $B_t\in\reals^{n\times m}$, and
$c_t\in\reals^n$ are random.  We assume the dynamics are time-invariant, \ie,
$(A_t, B_t, c_t)$ are independent and identically distributed (IID) for
different values of $t$. The initial state $x_0$ is also random, independent of
all $(A_t, B_t, c_t)$. When $A_t$, $B_t$, or $c_t$ are not random, \ie,
constant, we write them as $A$, $B$, or $c$.

\PAR{Certainty-equivalent dynamics}
We denote the expectations of the dynamics matrices as $\bar A = \Expect A_t$,
$\bar B=\Expect B_t$, and $\bar c = \Expect c_t$. We refer to the dynamical
system with the matrices replaced by their expectations,
\[
z_{t+1} = \bar Az_t + \bar B_t v_t + \bar c \quad t=0,1, \ldots,
\]
with initial condition $z_0=\Expect x_0$, as the certainty-equivalent system
(with state $z_t \in \reals^n$ and input $v_t \in \reals^m$).

\PAR{State-feedback policy}
We consider the time-invariant state feedback policy
\[
u_t = \phi(x_t), \quad t=0,1,\dots,
\]
where $\phi:\reals^n\to\reals^m$ is the policy that maps the state to the input.
The closed-loop system dynamics are
\[
x_{t+1} = A_t x_t + B_t\phi(x_t) +c_t, \quad t=0,1,\ldots,
\]
which defines a stochastic process for the state $x_t$.

\PAR{Stage cost}
The stage cost is a function $g:\reals^n\times\reals^m\to\reals\cup\{\infty\}$,
where $g(x_t,u_t)$ is the cost at time $t$. The stage cost $g$ imposes
constraints by taking infinite values at disallowed state-input pairs
$(x_t,u_t)$. We assume that the stage cost is a closed convex function. Note
that the cost function does not depend on time, \ie, it is time-invariant.

In some applications the cost is also random, \eg, of the form $\tilde
g_t(x_t,u_t)$, where $\tilde g_t$ is IID, and independent of
$A_t,B_t,c_t$, and therefore also of $x_t$.   Since we will work with the
expected value of the stage cost, we can handle this situation by taking
$g(x,t) = \Expect \tilde g_t (x,t)$, where the expectation is over the
random stage cost. For simplicity we assume that this expectation may be
computed analytically. In other cases, the expectation may be approximated, for
example using a sample average.

\PAR{Average cost}
The infinite-horizon average cost is given by
\BEQ\label{e-average-cost} J =
\lim_{T\to\infty} \frac{1}{T+1}\sum_{t=0}^{T} \Expect g(x_t,u_t).
\EEQ
Here, we assume that the limit and expectations exist.

We exclusively consider the average-cost problem, and do not consider the
closely-related discounted infinite horizon problem and finite horizon problem,
which may have time-varying stage cost. However, our approach is readily
extended to those problem settings, as discussed in \S\ref{s-extensions}.

\PAR{Convex stochastic control problem}
The convex stochastic control problem is to choose the policy $\phi$ so as to
minimize the cost $J$. We will denote an optimal policy as $\phi^\star$, and
assume that it exists. We let $J^\star$ denote the optimal value, \ie, the cost
$J$ with an optimal policy. The data in this problem are the distributions of
$(A_t,B_t,c_t)$ (which do not depend on $t$), the distribution of $x_0$, and the
stage cost function $g$.

\subsection{Dynamic programming}\label{s-dp}

The optimal control problem is readily solved, at least in principle, using
dynamic programming (DP)
\cite{Bellman1954,Pontryagin1987,BertsekasS1996,Puterman2014,Bertsekas2012}. An
optimal policy may be expressed in terms of a so-called Bellman or optimal value
function $V^\star : \reals^n \to \reals \cup \{\infty\}$, which roughly speaking
represents the optimal long-term cost of being in a given state.

An optimal policy can be expressed in terms of a value function as
\BEQ\label{e-opt-policy}
\phi^\star(x) = \argmin_u \left( g(x,u) + \Expect
V^\star(A_tx + B_tu + c_t) \right).
\EEQ
If there are multiple minima, we can arbitrarily choose one. The first term in
the quantity that is minimized is the immediate stage cost incurred by the input
choice $u$. The second term reflects the optimal expected long-term cost of
starting from the next state. The optimal policy balances these two costs.

The policy does not change when we add a constant to a value function. Without
loss of generality we can remove this ambiguity by insisting that
$V^\star(x^{\text{ref}}) = 0$, where $x^\text{ref}$ is a reference state (for
which there is an optimal value function with finite value). The value function 
\[
V^\text{rel}(x) = V^\star(x)-V^\star(x^\text{ref})
\]
is sometimes called a \emph{relative value function}.

\PAR{Bellman operator}
It can be shown that a value function $V^\star$ and the optimal cost $J^\star$
satisfy \BEQ\label{e-average-cost-bellman} V^\star + J^\star = \mathcal T
V^\star, \EEQ where $\mathcal T$ is the \emph{Bellman operator}, given by
\BEQ\label{e-bellman-op}
\left(\mathcal T h\right)(x) = \min_u \left( g(x,u)
+ \Expect h(A_tx + B_tu + c_t) \right),
\EEQ
for $h: \reals^n\to\reals\cup\{\infty\}$.

It follows that a relative value function is a fixed point of the Bellman
operator $\mathcal T$, \ie,
\[
V^\mathrm{rel} = \mathcal T V^\mathrm{rel}.
\]
This fixed point condition implies \eqref{e-average-cost-bellman}, with optimal
cost $J^\star = \mathcal TV^\star(x^\mathrm{ref})$.

\PAR{Value iteration}
The relative value function $V^\mathrm{rel}$ may be found by fixed point
iteration. Under certain technical conditions, the so-called value iteration (or
relative value iteration) \BEQ\label{e-value-iter} V^{k+1} = \mathcal T V^k -
\mathcal TV^k(x^\text{ref}), \\ \quad k=1,2,\dots \EEQ converges, \ie, $V^k -
V^k(x^\mathrm{ref}) \to V^\mathrm{rel}$ and $\mathcal T V^k(x^\text{ref}) \to
J^\star$ \cite{White1969,Puterman2014, Bertsekas2012}.

For future reference we mention a variation on value iteration called
\emph{damped value iteration}, which has the form
\BEQ\label{e-damped-value-iter} V^{k+1} = \rho_k \left(\mathcal T V^k - \mathcal
T V^k(x^\text{ref})\right) + (1-\rho_k) V^k, \quad k=1,2,\ldots,
\EEQ
where $\rho_k \in (0,1]$ with $\sum_k \rho_k(1-\rho_k) = \infty$. Damped value
iteration also satisfies $V^k - V^k(x^\mathrm{ref}) \to V^\mathrm{rel}$ and
$\mathcal T V^k(x^\text{ref}) \to J^\star$ under certain technical conditions.

\PAR{The value function is convex}
The Bellman operator \eqref{e-bellman-op} maps convex functions to convex
functions, since expectation and partial minimization preserve convexity (see,
\eg, \cite[\S3.2.1, \S3.2.5]{BoydV2004}). With any convex $V^1$ (\eg, the zero
function), it follows that all iterates of value iteration are convex, which
implies that its limit $V^\star$ is convex.

One implication is that evaluating the policy \eqref{e-opt-policy}, \ie,
minimizing 
\[
g(x,u) + \Expect V^\star(A_tx + B_tu + c_t)
\]
over $u$, is a convex optimization problem.  To see this, we observe that $A_t
x+ B_t u + c_t$ is an affine function of $u$, so by the affine pre-composition
rule, $V^\star(A_tx+b_tu+c_t)$ is a convex function of $u$. Adding 
this to $g(x,u)$ and taking expectation preserve convexity, so the function
that is minimized is a convex function of $u$.

Since evaluating the policy \eqref{e-opt-policy} involves solving a convex
optimization problem, we refer to it as a \emph{convex optimization control
policy}.


\PAR{Linear quadratic regulator}
The dynamic programming approach can only be carried out in practice in special
cases.  
The most widely known example is when the stage cost is a (convex) quadratic
function, in which case the optimal control problem is called the \emph{linear
quadratic regulator} (LQR). For LQR the Bellman operator preserves convex
quadratic functions, so it follows that the limit $V^\star$ is also convex
quadratic, and the optimal policy is affine, \ie, $\phi^\star(x) = Kx + l$,
where $K \in \reals^{m \times n}$ and $l \in \reals^m$ (see
\cite{BarrattB2021}). Value iteration for LQR can be carried out using basic
linear algebra operations, and so is tractable. Most importantly we have a
practical way to represent the Bellman iterates, and also their limit, by a
finite set of parameters, the coefficients of a quadratic function.

\PAR{Dynamic programing in the general case}
Beyond the special case of LQR described above, there are a handful of other
very specific stochastic control problems that are tractable to solve. These
cases follow the same general story line as LQR: There is a class of functions
that is preserved under the Bellman operator. One example is Merton's portfolio
problem, which considers the allocation of wealth between various assets over
time, and admits a closed-form solution \cite{Merton1969}. Problems with a
finite state space may, in principle, be solved by DP, by representing the value
function with a table of values. This is referred to as the tabular case
\cite{SuttonB2018}. When the state space is continuous but low-dimensional, say,
with $n\le 4$, the region of interest in the state space may be represented
using a finite number of points, for example a uniform grid. Tabular DP may then
be used, in combination with an interpolation over those points, to give a good
approximation of the value function. However, this approach does not scale to
problems with larger state dimension, since the number of points needed to
represent the value function to a given accuracy grows exponentially with the
state dimension.

The challenge in carrying out dynamic programming in more general cases is
simple: There is no practical way to represent an arbitrary convex function on
$\reals^n$.

\subsection{Certainty-equivalent steady-state optimal state-input pair}

For many stochastic control problems, certainty-equivalent approximations may be
used to obtain heuristic policies without dynamic programming. In this section
we explain the idea of an optimal steady-state certainty-equivalent optimal
state-input pair. We start by making two very crude approximations of the
stochastic control problem. First, we ignore all uncertainty by replacing the
dynamics matrices with their mean values (also called certainty-equivalent).
Second, we assume that the system is in steady-state, with constant state $z\in
\reals^n$ and constant input $v\in \reals^m$, \ie, $z=\bar A z + \bar B v + \bar
c$. Then we choose $z$ and $v$ to minimize the objective, which with the
assumptions above reduces to $g(z,v)$. Thus we solve the convex optimization
problem \BEQ\label{e-ce-ss-prob}
\begin{array}{ll} \mbox{minimize} & g(z,v) \\
\mbox{subject to} & z= \bar A z+ \bar B v + \bar c,
\end{array}
\EEQ with variables $z\in \reals^n$ and $v \in \reals^m$. We refer to a solution
of this problem $(z^\star,v^\star)$ as a certainty-equivalent steady-state
optimal (CE-SSO) state-input pair, and denote it as $(x^\mathrm{sso},
u^\mathrm{sso})$. For some problems, such as the example considered in
\S\ref{s-commitments}, the constant policy $\phi(x) = u^\mathrm{sso}$ is a
reasonable heuristic.

\subsection{Certainty-equivalent model predictive control}\label{s-ce-mpc}
Certainty-equivalent model predictive control (CE-MPC) is another heuristic
policy for stochastic control \cite{GarciaPM1989,BorrelliBM2017}. CE-MPC is not
our focus, but the methods of this paper can also be used to develop a good
CE-MPC policy.

To evaluate the CE-MPC policy $\phi^\mathrm{mpc}(x)$, we solve an $H$-step ahead
planning problem with certainty-equivalent dynamics. The planning problem is
\BEQ\label{e-mpc-prob}
\begin{array}{ll} \mbox{minimize} & \frac{1}{H+1} \sum_{\tau=1}^{H}
g(z_\tau,v_\tau) + V^\mathrm{mpc}(z_{H+1})\\
\mbox{subject to} & z_{\tau +1} = \bar A  z_\tau + \bar Bv_\tau + \bar c, \quad
\tau = 1, \ldots, H\\
& z_1=x,
\end{array}
\EEQ
with variables $z_1, \ldots, z_{H+1}$ and $v_1, \ldots, v_H$. The CE-MPC
policy is then $\phi^\mathrm{mpc} (x) = v_1^\star$, the first input of an
optimal trajectory of the MPC planning problem. \eqref{e-mpc-prob}.

In the CE-MPC problem \eqref{e-mpc-prob}, $V^\mathrm{mpc}$ is called the
terminal cost.  It can be chosen to be zero (particularly when $H$ is large
enough), or the indicator function of $x^\mathrm{sso}$, an optimal
certainty-equivalent steady-state state. Another very good choice is $\hat V$,
an approximation of the value function, which can be found by the methods of
this paper.

\section{Quadratic approximate dynamic programming}\label{s-qadp}

\subsection{Approximate dynamic programming}
In this paper, we consider ADP policies that replace the optimal value function
$V^\star$ in \eqref{e-opt-policy} with a convex approximation $\hat V$. The ADP
policy is of the form
\BEQ\label{e-adp-policy}
\hat\phi(x) = \argmin_u \left(g(x,u) + 
\Expect \hat V\left(A_t x + B_tu + c_t\right)\right).
\EEQ
(We omit the constant or offset term since it does not affect the associated
policy.)
If there are multiple minima, we can arbitrarily choose one. When $\hat V$ is a
convex quadratic function, we refer to \eqref{e-adp-policy} as a QADP policy.

ADP is a heuristic that addresses the issue mentioned above, that there is no
practical way to represent an arbitrary convex function on $\reals^n$
\cite{BellmanD1959,Munos2007,Bertsekas2012}. The approximate value function
$\hat V$ is chosen to approximate $V^\star$ in some sense, and to make
evaluating the policy \eqref{e-adp-policy} tractable. Evaluating $\hat\phi$ is
always a convex optimization problem; depending on the form of $g$ and $\hat
V$, the expectation can simplify and the problem can reduce to a common form,
such as a quadratic program (QP). When it is not possible to evaluate the
expectation in the policy exactly, we can use an estimate obtained by replacing
the expectation with a suitable sample average, \ie, a Monte Carlo approximation
\cite{Bertsekas2012}. ADP often works well in practice, even in cases when $\hat
V$ is not a particularly good approximation of $V^\star$
\cite{KeshavarzB2014,AgrawalBBS2020}.

\subsection{Quadratic approximate value functions}\label{s-quad-value-functions}

In this paper we focus exclusively on quadratic approximate value functions of
the form
\BEQ\label{e-quad-value-function}
\hat V(x) = \frac 1 2 \BBM x\\
1\EBM^T \BBM P & p \\ p^T & 0 \EBM \BBM x \\ 1 \EBM =
\frac 1 2 x^T Px + p^Tx,
\EEQ
where $P \succeq 0$, \ie, $P\in \symm_{+}^n$, the set of symmetric positive
semidefinite (PSD) $n \times n$ matrices.



The QADP policy associated with $\hat V$ is parametrized by the $n \times n$ PSD
matrix $P$ and $n$-vector $p$, which we collectively refer to as $\theta =
(P,p)$. All together, the parameter $\theta$ contains
\BEQ\label{e-qadp-parameter-count}
n(n+1)/2 + n = (1/2)n^2 + (3/2)n
\EEQ
scalar parameters, which has order $n^2$. We define $\Theta = \{ \theta \mid P
\succeq 0\}$, the set of parameters for which $\hat V$ is convex.

\subsection{Properties of QADP policies}\label{s-qadp-properties}

We now consider several properties of the QADP policies which will be useful in
the sequel.

\PAR{Simplifying the expectation}
The QADP policy can be simplified, since the expectation of a quadratic function
can be expressed analytically in terms of the first and second moments of its
argument. Thus we have
\BEQ\label{e-quad-expectation}
\Expect \hat V(A_tx + B_tu + c_t) = 
\frac 1 2 \BBM u \\ 1\EBM^T 
\BBM M & m \\ m^T & \mu(x) \EBM
\BBM u \\ 1 \EBM,
\EEQ
where
\[
\arraycolsep=1.5pt
\begin{array}{rl}
M &= \Expect B_t^TPB_t, \\[3pt]
m &= \Expect B_t^T(PA_tx+Pc_t) + \bar B^T p, \\[3pt]
\mu(x) &= x^T\Expect(A_t^TPA_t)x + 2x^T\Expect(A_t^TPc_t) + 2x^T\bar A^Tp
+ 2p^T\bar c + \Expect c_t^TPc_t.
\end{array}
\]

Note that $\mu(x)$ depends on $x$, and therefore is not constant, but the other
coefficients $M$ and $m$ are constant and depend only on the first and second
moments of $A$, $B$, $c$ (and $P$ and $p$). These formulas are derived in
\S\ref{a-expect-quad-fns}. Finally, we observe that $M$, $m$, and
$\mu(x)$ are linear functions of $\theta$.

\PAR{Evaluating the policy}
Since $g(x,u)$ is convex, evaluating the quadratic ADP policy reduces to
solving a deterministic convex optimization problem. When in addition
$g(x,u)$ is QP-representable, \ie, a convex quadratic function plus a convex
piecewise linear function, plus the indicator function of linear inequality and
equality constraints, evaluating the QADP policy reduces to solving a QP
\cite{WangB2010}.

\PAR{Gradient of the Bellman operator image}
Given convex quadratic $\hat V$, we may evaluate $\mathcal T\hat V(x)$, the
Bellman operator applied to $\hat V$ at any state $x$, by solving the convex
optimization problem associated with the QADP policy. We can also compute
$\nabla\mathcal T\hat V(x)$, where it is differentiable, and a subgradient
otherwise.

To do this, we represent $\mathcal T \hat V(x)$ as the optimal value of the
convex optimization problem
\BEQ\label{e-bellman-op-augmented}
\begin{array}{ll}
\mbox{minimize} & g(\tilde x,u) + \Expect \hat V(A_t \tilde x + B_tu + c_t) \\
\mbox{subject to} & \tilde x = x,
\end{array}
\EEQ
where we have introduced the variable $\tilde x$. Let $\nu^\star(x) \in\reals^n$
represent the optimal Lagrange multiplier associated with the constraint $\tilde
x = x$. Then, we have $\nabla\mathcal T \hat V (x) = -\nu^\star(x)$ when the
gradient exists \cite[\S 5.6]{BoydV2004}. Otherwise, $-\nu^\star(x)$ is a
subgradient, \ie, $-\nu^\star(x)\in\partial\mathcal T \hat V(x)$.

\section{Value-gradient iteration}\label{s-qfvgi}

\subsection{Fitted value iteration}\label{s-fvi}
We begin by reviewing fitted (or projected) value iteration (FVI), which is an
approximation of value iteration \cite{BellmanD1959, KeshavarzB2014,
Bertsekas2012}. The issue with value iteration is that in practice, we cannot
exactly represent the function $\mathcal T V^k$ in the update
\eqref{e-damped-value-iter}. FVI addresses this by restricting all approximate
value function iterates $V^k$ to be convex quadratic functions.

In the $k$th iteration, we choose a set of states $x^1,\ldots,x^N$, and evaluate
$\mathcal T V^k(x^i)$ for each $i=1,\ldots,N$. We can evaluate each $\mathcal T
V^k(x^i)$ by evaluating \eqref{e-bellman-op}, which is a convex optimization
problem. Then, we fit a convex quadratic function $V^{k+1/2}$ to those points,
such that
\[
V^{k+1/2}(x^i) \approx \mathcal T V^k(x^i), \quad i=1,\ldots,N.
\]
This leads to the damped fitted value iteration update
\BEQ\label{e-fitted-value-iter}
V^{k+1} = \rho_k V^{k+1/2} + (1-\rho_k) V^k, \quad k=0,1,\ldots,
\EEQ
which generates a sequence of convex quadratic functions $V^k$, with associated
QADP policies.

\PAR{Fitting convex quadratic functions}
One method for finding parameters $\theta=(P,p)$ for the convex quadratic
function $V^{k+1/2}$ is to fit it to a set of points. We first evaluate
$v^i = \mathcal T V^k(x^i)$ for each $i=1,\ldots,N$, and then solve the fitting
problem
\BEQ\label{e-fvi-fitting}
\begin{array}{ll}
\mbox{minimize} & (1/N)\sum_{i=1}^N
L\left(V^{k+1/2}(x^i) + c - v^i\right) + r(\theta) \\
\mbox{subject to} & \theta\in\Theta,
\end{array}
\EEQ
with variables $\theta$ and $c\in\reals$, where $c$ is a scalar offset. Here $L
: \reals\to\reals$ is a convex fitting loss function, and $r :
\symm^n\times\reals^n\to\reals\cup\{\infty\}$ is a convex regularization
function, with infinite values used to impose (convex) constraints on $\theta$.
This is a convex optimization problem, since $V^{k+1/2}(x^i)$ is a linear
function of $\theta$. Possible choices for $L$ include the squared loss or the
robust Huber loss \cite{Huber1992}, given by
\BEQ\label{e-huber-scalar}
L^{\textrm{hub}}(z) = \begin{cases}
(1/2) z^2 & |z| \leq M \\
M(|z| - M/2) & |z| > M.
\end{cases}
\EEQ
The Huber loss is a more robust alternative to the square loss, in the presence
of outliers. Possible choices for $r$ include $\ell_2$ regularization and prior
knowledge constraints, and are discussed in \S\ref{s-reg}. For simplicity, we
consider the standard Huber function, which transitions from the quadratic to
absolute value at $M=1$. In general, $M$ may be tuned by cross-validation, using
a procedure similar to that described in \S\ref{s-reg}.

\PAR{Convergence}
Convergence guarantees for FVI are available when the approximation error of
$\nabla\mathcal T V^k$ is small enough \cite{Munos2007,Bertsekas2012}. However,
unlike value iteration, FVI is not guaranteed to converge in general
\cite{Baird1995,TsitsiklisVR1996}. Nevertheless, with an appropriate
approximation $\hat\nabla\mathcal T V^k$ and damping parameters $\rho_k$, FVI
can often find policies with good performance in practice.

\subsection{Value-gradient iteration}
VGI is a special case of FVI, where we fit $V^{k+1/2}$ using gradients instead
of values. In \S\ref{s-qadp-properties}, we showed that we can evaluate $\nabla
\mathcal T V^k(x)$ at any state $x$ where $\mathcal T V^k$ is differentiable, by
evaluating a particular optimal Lagrange multiplier. Therefore, we can find
$V^{k+1/2}$ by fitting its gradient.

That is, we choose $V^{k+1/2}(x) = (1/2)x^TP x + p^T x$ such that $P\succeq 0$
and
\[
\nabla V^{k+1/2}(x) = Px^i + p \approx \nabla\mathcal T V^k(x^i),
\quad i=1,\ldots,N.
\]
Once we have found $V^{k+1/2}$, we apply the damped update
\eqref{e-fitted-value-iter} to generate the next iterate $V^{k+1}$. Like in
standard FVI, this generates a sequence of convex quadratic functions $V^k$,
with associated QADP policies.

\PAR{Fitting the gradient}
In this case, we fit an affine function $\nabla V^{k+1/2}(x) = Px + p$ to a set
of points, subject to the constraint that $P$ is symmetric positive
semidefinite. In each iteration, we evaluate $g^i = \nabla \mathcal T V^k(x^i)$
for each $i=1,\ldots,N$, and then solve the fitting problem
\BEQ\label{e-vgi-fitting}
\begin{array}{ll}
\mbox{minimize} & (1/N)\sum_{i=1}^N
L\left(\nabla V^{k+1/2}(x^i) - g^i\right) + r(\theta) \\
\mbox{subject to} & \theta\in\Theta,
\end{array}
\EEQ
with variables $\theta$. Here $L : \reals^n\to\reals$ is a multivariate convex
fitting loss function, and $r$ is, like in \eqref{e-fvi-fitting}, a convex
regularization function. This is also a convex optimization problem, since
$\nabla V^{k+1/2}(x^i)$ is a linear function of $\theta$.

Possible choices for $L$ include the squared $\ell_2$ norm and the circular
Huber loss
\BEQ\label{e-huber-loss} L^{\textrm{hub}}(z) = 
\begin{cases} (1/2) \|z\|_2^2 &
\|z\|_2 \leq M \\
M(\|z\|_2 - M/2) & \|z\|_2 > M,
\end{cases}
\EEQ
which extends the scalar Huber loss \eqref{e-huber-scalar} to the multivariate
case. Like in the scalar case, the circular Huber loss is a more robust
alternative to the square function, in the presence of outliers.

\PAR{Choice of sampling points}
An important consideration is the choice of the state samples values $x^1,
\ldots, x^N$ at which we evaluate the policy and $\mathcal T V^k(x^i)$. Ideally
the samples should reflect the states that the system is likely to be in, \ie,
samples from the steady-state distribution of $x_t$ under the policy $\phi^k$.

To accomplish this we choose the sample points by simulating the current policy
for $N$ steps, using the current policy $\phi^k$. In the first iteration $k=1$,  
we initialize the simulation at a state chosen at random. In subsequent
iterations, we initialize the simulation at the last state in the previous
iteration.

\subsection{Regularization, constraints, and lower bounds}\label{s-reg}

Prior information, if available, can be incorporated as regularization terms or
constraints in the fitting problem, through the function $r(\theta)$ in the
fitting problem \eqref{e-vgi-fitting}. Constraints and lower bounds may be
imposed by setting $r$ to have value $\infty$ when $\theta$ is not consistent
with the prior information. We now describe a nonexhaustive list of
possibilities that may be combined to form $r(\theta)$.

\PAR{Ridge regularization}
We may add an $\ell_2$ penalty on the parameters of the value function
\[
r(\theta) = \lambda\left(\|P\|_F^2 + \|p\|_2^2\right),
\]
where $\lambda > 0$ is a scalar regularization parameter and $\|\cdot\|_F$
denotes the Frobenius norm. The $\ell_2$ regularization ensures that the fitting
problem is well-posed and helps mitigate overfitting, and is sometimes referred
to as Tikhonov or ridge regularization \cite{TikhonovA1977,HastieTF2009}.

The parameter $\lambda$ is typically chosen using use out-of-sample or
cross-validation. To do this we divide the fitting data $(x^i,v^i)$ into two
sets, the training data and the validation data. We fit $V$ using the training
data, for a range of values of $\lambda$, typically on a log scale with upper
limits $\lambda^\mathrm{max}$ and $\lambda^\mathrm{max}$, and then evaluate the
average loss on the validation data for each value of $\lambda$. We then choose
a value that gives near minimum validation error, with a preference for larger
values, \ie, more regularization. This approach is often referred to as grid
search. A more thorough method is to use cross-validation \cite{HastieTF2009},
and more sophisticated search methods for evaluating scaling parameters may also
be considered; see, for example, \cite{JamiesonT2016}.

\PAR{LASSO regularization}
The $\ell_1$ penalty
\[
r(\theta) = \lambda\left(\sum_{i,j=1}^n |P_{ij}| + \|p\|_1\right)
\]
with regularization parameter $\lambda > 0$ is known as LASSO
\cite{HastieTF2009}. This regularization is similar to ridge regression in that
both shrink the values of the parameters; however, the LASSO is more likely to
produce sparse solutions, \ie, $P$ and $p$ with zero-valued entries. Therefore,
the LASSO regularization can be particularly useful for weakly coupled systems.

Like with ridge regression, the value of $\lambda$ may be tuned using
out-of-sample or cross-validation. When multiple regularization terms are used,
we can use the same strategy to find a good set of values for each
regularization parameter. For example, the case where both ridge and LASSO
regularization are employed is known as the elastic net \cite{ZouH2005}. In this
case, the aforementioned grid search strategy may be used to select the two
regularization parameters jointly.

\PAR{Symmetry}
In some cases, we may know that the value function $V$ should be symmetric, \ie,
$V(x) = V(-x)$ for any $x\in\reals^n$. The LQR example considered in
\S\ref{s-box-lqr}, for example, satisfies this property. For quadratic
approximate value functions, symmetry may be implemented by the constraint
$p=0$.

\PAR{Fixed minimizer}
When we can identify a point $x^\star$ in the state space that seems to be the
best, we may include the constraint $\argmin_x V(x) = x^\star$ to the fitting
problem. This is equivalent to the linear equality constraint $Px^\star + p =
0$. A special case is when $V$ is constrained to be symmetric, in which case
$V(x)$ is minimized at zero.

\PAR{Lower bounds}
In some cases, a quadratic pointwise lower bound
\[
V^{\mathrm{lb}}(x) = 
\frac 1 2 x^T P^{\mathrm{lb}}x + (p^{\mathrm{lb}})^T x
= \frac 1 2
\BBM x \\ 1 \EBM^T
\BBM P^{\mathrm{lb}} & p^{\mathrm{lb}} \\
(p^{\mathrm{lb}})^T & 0
\EBM
\BBM x \\ 1 \EBM
\]
on $V^\star$ is available up to an additive constant, and may be included as an
additional constraint. This may be done by introducing an additional variable
$s$, and imposing the pointwise constraint $V + s \geq V^\mathrm{lb}$. This can
be expressed as the convex constraint
\BEQ\label{e-quad-lb}
\BBM P - P^{\text{lb}} & p - p^{\text{lb}} \cr
(p - p^{\text{lb}})^T & s
\EBM \succeq 0,
\EEQ
as shown in \S\ref{a-quad-lower-bound}. Since $P^\mathrm{lb} \succeq 0$, this
constraint implies that $P \succeq 0$. So when we add a quadratic lower bound
constraint to the fitting problem, we no longer need the constraint $P \succeq
0$.

In many cases we can form a convex quadratic lower bound $V^{\text{lb}}$ on the
true value function $V^\star$. In the simplest case we can take $V^\star=0$ when
the stage cost is nonnegative. Another method is to form an LQR relaxation of
the problem, \ie, to replace $g$ with a quadratic lower bound, for example,
by ignoring constraints on $u$. The resulting LQR problem can be solved exactly,
and its value function $V^\mathrm{lqr}$ is a lower bound on $V^\star$. More
sophisticated methods for computing a lower bound on the value function involve
solving a convex optimization problem \cite{WangB2009} or a series of convex
problems \cite{OWB2011}.

When the dynamics matrices $A_t$ and $B_t$ are random, a simpler lower bound may
be found by considering the (deterministic) LQR relaxation of the CE problem;
see \S\ref{a-ce-lb}. 

\PAR{Policy interpolation}
Suppose we have a set of states $x^1,\ldots, x^B$, and require that the policy
takes on corresponding values $u^1,\ldots,u^B$, \ie,
\[
\phi(x^j) = u^j, \quad j=1,\ldots,B.
\]
This condition may be written as
\BEQ\label{e-policy-interp}
0 \in \partial g(x^j, u^j) + \nabla \Expect V^k(A_t x^j + B_t u^j + c_t),
\EEQ
where $\partial g(x^j, u^j)$ is the set of subgradients of $g(x,u)$ with
respect to $u$, evaluated at $(x^j, u^j)$.

In some cases, this constraint has a simple representation. For example, if the
stage cost may be written in the form
\[
g(x,u) = h(x, u) + I\left((x, u) \in C\right)
\]
where $h$ is differentiable and $I\left((x, u) \in C\right)$ is the indicator
function of a polyhedral set $C$, then the constraint may be written as a linear
inequality constraint on the parameters $P$ and $p$. First, note that
\[
\partial g(x^j, u^j) = 
\nabla h(x^j, u^j) + \partial I\left((x, u) \in C\right),
\]
where $\partial I\left((x, u) \in C\right)$ is the normal cone to $C$ at $(x^j,
u^j)$. Since $C$ is a polyhedron the normal cone is also a polyhedron \cite[\S
23]{Rockafellar1970}, \ie, representable by a set of linear inequality
constraints. Next, from \eqref{e-quad-expectation} we have
\[
\nabla \Expect V^k(A_t x^j + B_t u^j + c_t)
= \Expect(B_t^T P B_t)u^j + \Expect B_t^T(PA_t x^j + Pc_t) + \bar B^T p,
\]
which is a linear function of $P$ and $p$. Therefore, the policy
interpolation constraints  \eqref{e-policy-interp}
may be represented by a set of linear inequality constraints on $P$ and $p$.

\section{Extensions and variations} \label{s-extensions}

\subsection{Input-affine dynamics}
The methods presented in this paper can also be applied in cases where the
dynamics are nonlinear but input-affine. That is, the dynamics may be written in
the form
\[
x_{t+1} = f_t(x_t) + B_t(x_t)u_t,
\]
where $f_t : \reals^n \to \reals^{n}$ and $B_t : \reals^n \to \reals^{n\times
m}$ are random functions. We again assume that $(f_t,g_t)$ are IID for different
values of $t$. The affine dynamics described in \S\ref{s-cvx-stoch-control} are
a special case, where $f_t(x) = A_tx + c_t$ and $g_t(x) = B_t$.

In the input-affine case, the ADP policy \eqref{e-adp-policy} is of the form
\[
\hat\phi(x) = \argmin_u \left(g(x,u) +
\Expect \hat V\left(f_t(x) + g_t(x)u\right)\right).
\]
Since the dynamics are affine in $u$, the expected value $\Expect \hat V(f_t(x)
+ g_t(x)u)$ is also affine in $u$, when $\hat V$ is convex. When $\hat V$ is a
convex quadratic function of the form \eqref{e-quad-value-function}, the
expected value may be computed exactly, in terms of the first and second moments
of $f_t(x)$ and $g_t(x)$ \cite{KeshavarzB2014}. Hence, the policy can still be
evaluated by solving a convex optimization problem, and VGI can still be
performed in a similar manner.

\subsection{Alternative cost functions}\label{s-alternative-costs}
\PAR{Discounted infinite-horizon problem}
The mean discounted infinite-horizon cost is given by
\[
J =\sum_{t=0}^\infty \gamma^t \Expect g_t(x_t, u_t),
\]
where $\gamma\in(0,1)$ is a discount factor, and the sum and expectations are
assumed to exist. In this case, the value function $V^\star$ represents the
optimal cost-to-go, and the optimal policy is of the form
\[
\phi^\star(x) = \argmin_u \left(g(x,u) +
\gamma \Expect V^\star(A_t x + B_tu + c_t)\right).
\]
For the discounted infinite-horizon problem, VGI proceeds in the same way,
except with the Bellman operator defined as
\[
\left(\mathcal T h\right)(x) = \min_u \left( g(x,u)
+ \gamma\Expect h(A_tx + B_tu + c_t) \right),
\]
for $h: \reals^n\to\reals\cup\{\infty\}$.

\PAR{Finite-horizon problem}
In the finite-horizon problem, the cost is given by
\[
J = \sum_{t=0}^T \Expect g_t(x_t, u_t), 
\]
where the stage cost may be time-varying, and the expectations are assumed to
exist. In this case, the value function $V_t^\star$ depends on time, and may be
found using a backward recursion. The value iteration starts with
\[
V^\star_T(x) = \min_u g_T(x,u),
\]
and then proceeds as
\[
V^\star_t(x) = \mathcal T_t V^\star_{t+1}(x), \quad t=T, T-1,\ldots,0,
\]
where the Bellman operator at time $t$ is defined as
\[
\left(\mathcal T_t h\right)(x) = \min_u \left( g_t(x,u)
+ \gamma\Expect h(A_tx + B_tu + c_t) \right),
\]
for $h: \reals^n\to\reals\cup\{\infty\}$.

VGI proceeds similarly for the finite-horizon problem, using an analogous
function fitting approximation of the Bellman operator.

\subsection{Parallel simulations}\label{s-parallel-sim}

In VGI (and FVI in general), we select $N$ sample points by simulating the
current policy.  We can also select points from more
than one simulated trajectory. To do this we choose the sample points by
simulating $K$ different trajectories for $T$ steps each, using the current policy. In
iteration $k$, each of these $K$ trajectories gives us $T$ states at which we
evaluate the policy $\phi^k$, so all together we have $N=TK$ states and
associated evaluations of $\nabla \mathcal TV^k$ to use in the fitting problem
\eqref{e-vgi-fitting}.
One advantage of this method is that the $K$ trajectories can be evaluated
in parallel.

\section{Numerical examples}\label{s-examples}

In this section, we present three numerical examples, which involve a
box-constrained LQR problem, a commitment planning problem with an alternative
investments fund, and a supply chain optimization problem. Comparisons with
other ADP methods are given in \S\ref{s-comparisons}.

The code for the examples is available at \url{https://github.com/cvxgrp/vgi}.
The ADP policies and VGI method are implemented using CVXPY
\cite{DiamondB2016,AgrawalVDB2018}. In addition, the code generation tool
CVXPYgen \cite{SchallerBDASB2022} was used to create custom solvers for the ADP
policies, implemented in C. The experiments were performed on two cores of an
Intel Xeon E5-2640 CPU.

\subsection{Box-constrained linear quadratic regulator}\label{s-box-lqr}

We first consider a traditional linear quadratic regulator (LQR) problem. The
dynamics are time-invariant, and given by
\[
x_{t+1} = Ax_t + Bu_t + c_t,
\]
where $A\in\reals^{n\times n}$ and $B^{n\times m}$ are known and fixed, and
$c_t$ is an IID random variable with zero mean and covariance $\Expect c_tc_t^T
= C$. The stage cost is given by
\[
g(x,u) = x^TQx + u^TRu + I(-u^{\max} \le u \le u^{\max}),
\]
where $Q\succeq 0$, $R\succ 0$, and $u^{\max}>0$ is a maximum input magnitude,
in any component of the input.

For this problem, a lower bound $J^{\textrm{lb}}$ on the optimal cost and a
quadratic lower bound $V^{\textrm{lb}}$ on the optimal value function can be
found by solving a semidefinite program (SDP) \cite{WangB2009}. An upper bound
on the optimal cost may be found by evaluating the ADP policy using
$V^{\textrm{lb}}$ as the approximate value function.

\PAR{Numerical example} We consider a problem instance with $n=12$ and
$m=3$. The entries of $A$ are chosen IID from a uniform distribution on
$[-1,1]$. The matrix $A$ was then rescaled to have a maximum eigenvalue of 1.
The entries of $B$ are chosen IID from a uniform distribution on $[-0.5,0.5]$.
The process noise $c_t$ is normally distributed, with zero mean and covariance
$0.4I$. The stage cost parameters are given by $Q=I$ and $R=I$, and the maximum
input magnitude is $u^{\max}=0.4$.

\PAR{Results}
We carried out VGI for $40$ iterations, starting from the initial value function
$V^1(x) = x^TQx$. We included the symmetry constraint $p=0$ in the fitting step.
In each iteration, the fitting step was performed using $N=50$ fitting points,
obtained by simulating the current policy. The damping coefficient was fixed to
$\rho_k=0.5$.

\begin{figure}
\centering
\includegraphics[scale=0.6]{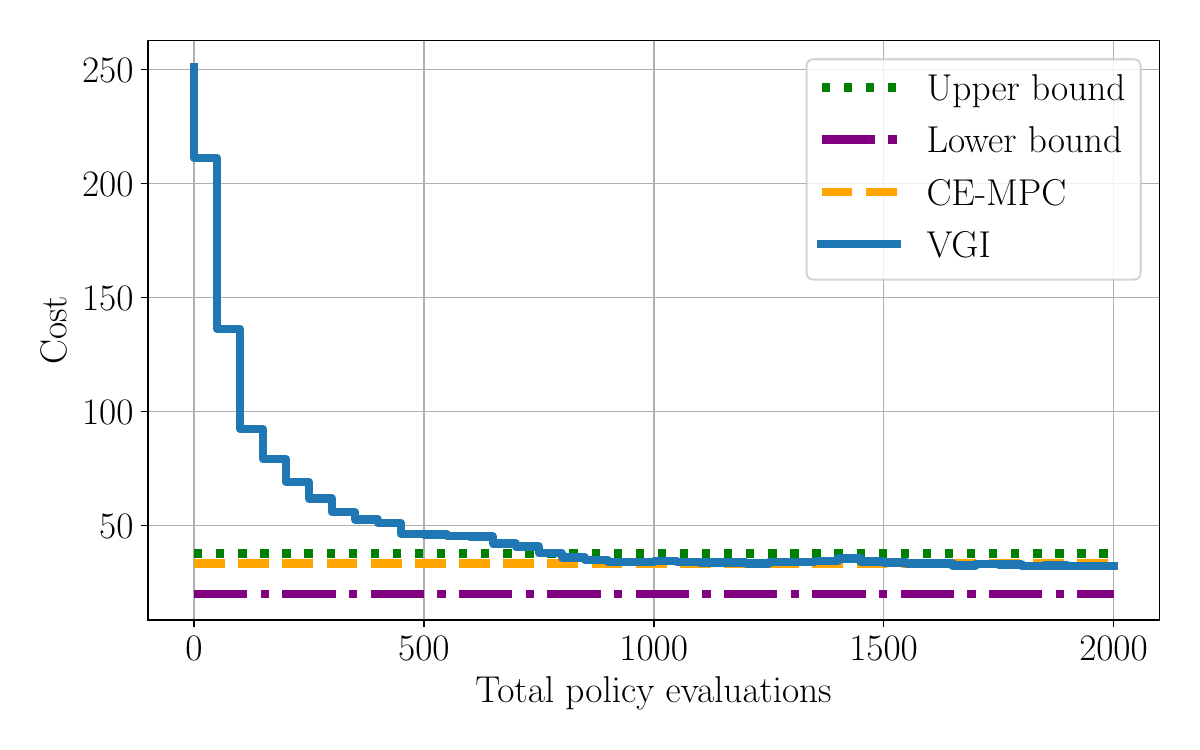}
\caption{VGI for the box-constrained LQR problem.}
\label{f-lqr-vgi}
\end{figure}

Figure \ref{f-lqr-vgi} shows the average cost versus the number of policy
evaluations used to generate the data for the fitting step. Also plotted are the
SDP-based upper and lower bounds \cite{WangB2009} and the average cost of the
CE-MPC policy with a horizon of $H=30$. In this example, VGI converges to a
slightly better cost than that of the CE-MPC policy.

\subsection{Commitments in an alternative investments fund}\label{s-commitments}
Our next example is a practical example, and more specific. We consider a fund
that invests in $m$ so-called alternative investment classes, such as venture
capital, infrastructure projects, direct lending, or private equity. Alternative
investments are found in the portfolios of insurance companies, retirement
funds, and university endowments. For more details, see
\cite{LuxenbergBvBCK2022} and the papers cited therein.

In each time period (typically quarters) $t=1,2, \ldots$, we make nonnegative
commitments to the $m$ alternative asset classes.  These are amounts we promise
to invest, in response to capital calls. Over the next few years, we put money
into the investments in response to capital calls, up to the amount of previous
commitments. We receive money from each the investments in later years through
distributions.  Neither the timing nor amounts of the capital calls and
distributions are directly under our control, except that the total of the
capital calls cannot exceed our total commitments for each asset class.

We first describe some critical quantities.
\begin{itemize}
\item $u_t \in \reals_+^m$ denotes the amounts that the investor commits in
period $t$, to each of the $m$ asset classes. (These commitments will be the
input in our stochastic control problem.)
\item $p_t \in \reals_+^m $ denotes the amounts that the investor pays in to the
investment in response to capital calls in period $t$.
\item $d_t\in \reals_+^m$ denotes the amount that the investor receives in
distributions from the investments in period $t$.
\item $n_t\in \reals_+^m$ denotes the net asset values (NAVs) of the investments
in period $t$.
\item $l_t\in \reals_+^m$ denotes the total amount of uncalled commitments, \ie,
the difference between the total so far committed and the total so far that has
been called. (This is a liability, so we use the symbol $l$.)
\end{itemize}
The units for all of these is typically millions of USD.

A simple dynamical model relating these variables is
\[
n_{t+1} = \diag(r_t) n_t + p_t - d_t, \quad 
l_{t+1} = l_t - p_t + u_t, \quad 
t=1,2 \ldots,
\]
where $r_t \in \reals_{++}^K$ is the vector of per-period total returns for the
asset classes, assumed to be IID with some known distribution such as
log-normal. In words: the value of each investment class in each period is
multiplied by its (random) return, increased by the amount paid in, and
decreased by the amount distributed; the total uncalled commitments is decreased
by the capital calls, and increased by new commitments. The calls and
distributions are modeled as
\[
p_t = \diag(\gamma^\text{call}_t) l_t, \quad
d_t = \diag(\gamma^\text{dist}_t)\diag(r_t) n_t, \quad t=1,2, \ldots,
\]
where $\gamma^\text{call}_t$ and $\gamma^\text{dist}_t$ are random variables in
$(0,1)^m$, called the call and distribution intensities. We will assume that
these are IID, and independent of $r_t$. In words: In each period and for each
asset class, a random fraction of the total liability is called, and a random
fraction of the NAV is distributed.

We can express the dynamics as a random linear dynamical system with state $x_t
= (n_t,l_t)\in \reals^{2m}$ and input $u_t \in \reals^m$, with dynamics matrices
\[
A_t = \left[ \begin{array}{cc}
\diag(r_t)(I-\diag(\gamma^\mathrm{dist}_t)) & \diag(\gamma^\mathrm{call}_t) \\
0 & I - \diag(\gamma^\mathrm{call}_t )\\
\end{array}\right], \qquad
B_t = \left[ \begin{array}{c} 0 \\ I 
\end{array}\right], \qquad c_t = 0.
\]

The goal is to choose commitments so as to reach and maintain a target asset
allocation $n^{\mathrm{tar}} \in\reals_+^m$, while penalizing deviations of the
commitments $u_t$ from the CE-SSO commitment $u^{\mathrm{sso}}\in\reals^m_+$. We
consider stage cost
\[
g(x_t, u_t) = \|n_t - n^{\mathrm{tar}}\|^2
+ \lambda \|u_t - u^{\mathrm{sso}}\|^2 + I(0\le u_t \le u^{\mathrm{max}}),
\]
where $\lambda > 0$ is a penalty coefficient and $u^{\mathrm{max}}\in\reals^m_+$
are the maximum allowable commitments to each of the asset classes. We take the
fixed input $u^{\mathrm{sso}}$ is a solution to the certainty-equivalent
steady-state problem \eqref{e-ce-ss-prob}, with the input cost term
$\lambda\|u_t - u^{\mathrm{sso}}\|^2$ removed from the stage cost.

For this problem, we find a quadratic lower bound $V^{\textrm{lb}}$ on the value
function by relaxing the constraints on the input $u_t$, replacing $A_t$ with
$\bar A$, and solving the certainty equivalent LQR problem.



\PAR{Numerical example}
We consider an example with $m=6$ asset classes. The returns $r_t$ are
distributed according to a log-normal distribution, \ie, $r_t = \exp(z_t)$, with
$z_t \sim \mathcal N(\mu,\Sigma)$. The parameters $\mu$ and $\Sigma$ were chosen
such that the mean quarterly returns have means
\[
\Expect r_t = \left(1.0, 1.1, 1.1, 1.0, 1.1, 1.1\right)
\]
and standard deviations
\[
\sigma_t = \left(0.1, 0.2, 0.2, 0.1, 0.2, 0.1 \right).
\]
This leads to annualized returns with means around $20\%$ and standard
deviations around $30\%$. The returns are correlated, with correlation matrix
\[
\mathrm{corr}(r_t) =
\left[\begin{array}{cccccc}
1 & -0.06 & -0.05 & 0.62 & -0.32 & -0.44 \\
-0.06 & 1 & -0.21 & 0.18 & 0.80 & -0.12 \\
-0.05 & -0.21 & 1 & 0.35 & -0.27 & -0.19 \\ 
0.62 & 0.18 & 0.35 & 1 & 0.18 & -0.15 \\
-0.32 & 0.80 & -0.27 & 0.18 & 1 & 0.37 \\
-0.44 & -0.12 & -0.19 & -0.15 & 0.37 & 1 \\
\end{array}\right].    
\]

The components of $\gamma^{\text{call}}_t$ and $\gamma^{\text{dist}}_t$ are
independent and beta-distributed, such that
$(\gamma_t^{\textrm{call}})_i\sim\mathrm{Beta}(\alpha^{\text{call}}_i,
\beta^{\text{call}}_i)$, where $\alpha^{\mathrm{call}}_i = 2$ for
$i=1,\ldots,m$, and
\[
\beta^{\mathrm{call}} = 
\left(10.3, 10.0, 12.9, 10.5, 11.8, 10.5\right).
\]
The distribution intensities were also beta distributed, such that
$(\gamma_t^{\textrm{dist}})_i\sim\mathrm{Beta}(\alpha^{\text{dist}}_i,
\beta^{\text{dist}}_i)$, where $\alpha^{\mathrm{dist}}_i = 3$ for
$i=1,\ldots,m$, and
\[
\beta^{\mathrm{dist}} =
\left(13.0, 12.7, 15.9, 12.8, 13.2, 14.2\right).
\]
These parameters lead to typical values of call and distribution intensities
around $0.14$ and $0.16$ respectively. The target asset values
$n^{\textrm{tar}}$ are chosen to be between 4 and 5, the maximum commitment is
$u^{\max} = 3$, and the penalty coefficient was $\lambda = 0.01$.

\begin{figure}
\centering
\includegraphics[scale=0.6]{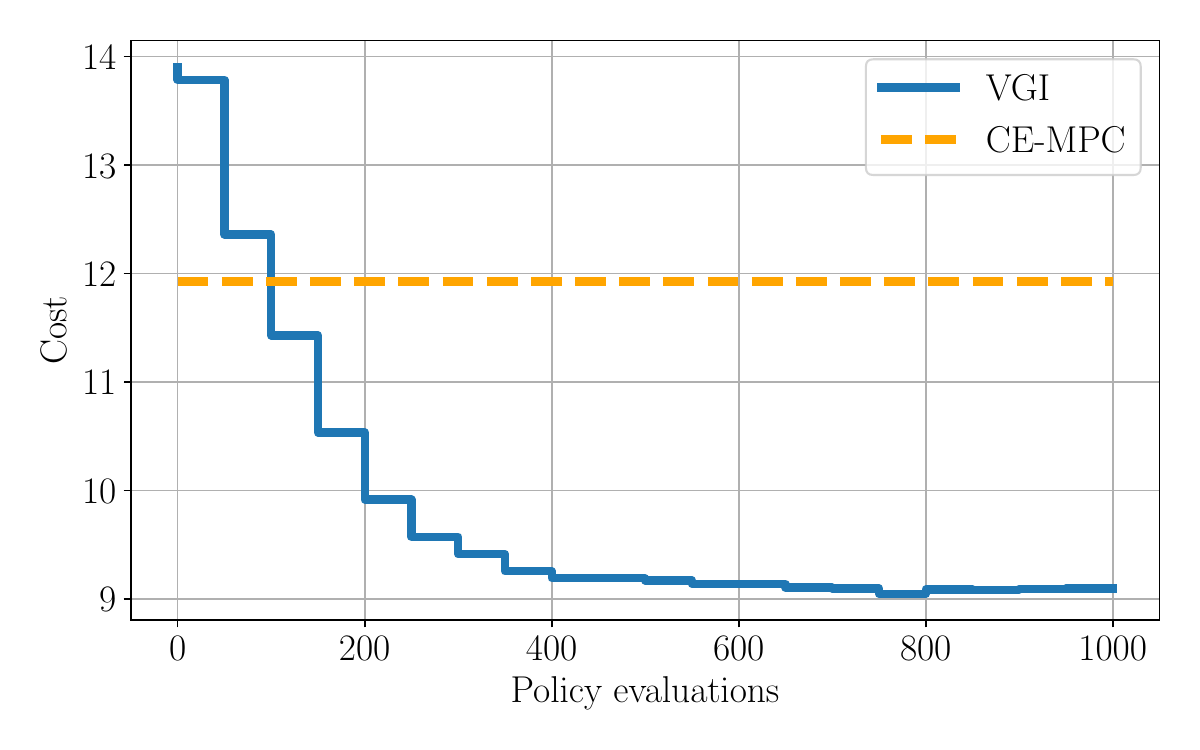}
\caption{VGI for the commitments planning problem.}
\label{f-commitments-vgi}
\end{figure}

\begin{figure}
\centering
\includegraphics[scale=0.6]{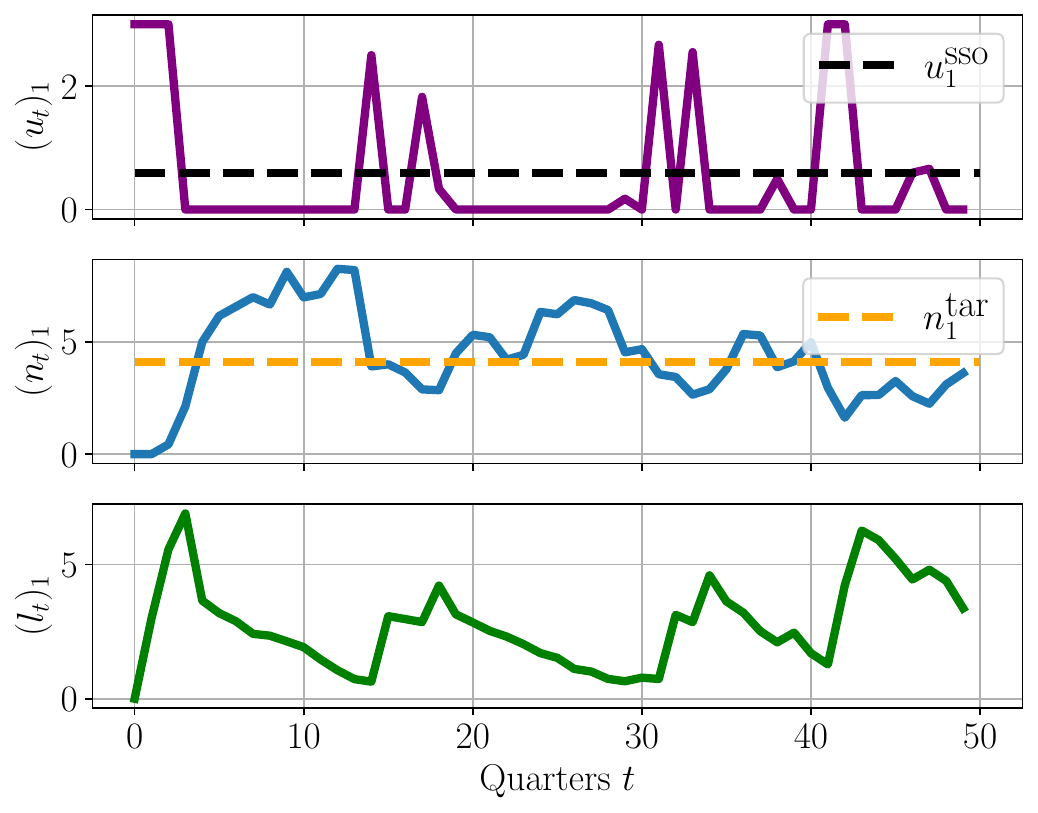}
\caption{Commitments, NAV, and uncalled commitments for one asset class. The
policy found by VGI makes commitments when the NAV dips below the target value.}
\label{f-commitments-trajectory}
\end{figure}
    
\PAR{Results}
We carried out VGI for $20$ iterations, starting from $V^1=V^{\mathrm{lb}}$.
In each iteration, the fitting step was performed using $N=50$ fitting points,
obtained by simulating the current policy. The damping coefficient was fixed to
$\rho_k=0.5$.

Figure \ref{f-commitments-vgi} plots the average cost versus the number of
policy evaluations used, along with the average cost of the CE-MPC policy with a
horizon of $H=30$. Our method converges to a policy that is 25\% better than the
CE-MPC policy. It is able to significantly outperform the CE-MPC policy because
it accounts for the correlation between the returns $r_t$. The CE-MPC policy, on
the other hand, only accounts for the average returns. The average costs were
computed by simulating the system for ten thousand steps.

Figure \ref{f-commitments-trajectory} shows an example trajectory of asset
value, liability, and commitments made for one of the six asset classes, using
the ADP policy found by VGI. The policy makes commitments when the asset value
dips below the target value.

\begin{figure}
\centering
\begin{tikzpicture}[
roundnode/.style={
circle, draw=blue!90, fill=blue!7.5, very thick, minimum size=10.5mm
},
squarednode/.style={
rectangle, draw=red!60, fill=red!5, very thick, minimum size=5mm
},
align=center,
node distance=1.5cm
]
\node[squarednode]      (source1)   {$\text{supplier 1}$};
\node[squarednode]      (source2)[below=of source1] {$\text{supplier 2}$};

\node[roundnode]        (node1)             [right=of source1]  {$h_1$};
\node[roundnode]        (node2)             [right=of node1]    {$h_3$};
\node[roundnode]        (node3)             [right=of source2]  {$h_2$};
\node[roundnode]        (node4)             [right=of node3]    {$h_4$};

\node[squarednode] (sink1) [right=of node2] {$\text{consumer 1}$};
\node[squarednode] (sink2)  [right=of node4] {$\text{consumer 2}$};

\draw[-{stealth[scale width=2]}] (source1.east) -- (node1.west);
\draw[-{stealth[scale width=2]}] (source2.east) -- (node3.west);

\draw[-{stealth[scale width=2]}] (node1.east) -- (node2.west);
\draw[-{stealth[scale width=2]}] (node3.east) -- (node4.west);

\draw[-{stealth[scale width=2]}] (node1.south east) -- (node4.north west);
\draw[-{stealth[scale width=2]}] (node4.north) -- (node2.south);

\draw[-{stealth[scale width=2]}] (node2.east) -- (sink1.west);
\draw[-{stealth[scale width=2]}] (node4.east) -- (sink2.west);
\end{tikzpicture}
    
\caption{Supply chain network.}
\label{f-sc-network}
\end{figure}
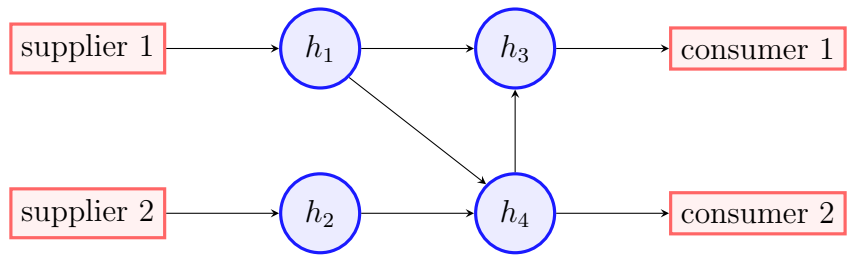
    
\subsection{Supply chain optimization}\label{s-sco}

In our final example, we consider the problem of shipping goods efficiently
across a network of warehouses to maximize profit. We consider a single-good,
multi-echelon supply chain with $\tilde n$ interconnected warehouses, which are
represented by nodes in a graph. There are $m$ directed links over which goods
can flow; $n_s$ links connect suppliers to nodes, $n_c$ links connect nodes to
consumers, and $m-n_s-n_c$ links connect nodes to each other.

The amount of good held at each node at time $t$ is represented by
$h_t\in\reals_+^{\tilde n}$. The prices at which we can buy the good from the
suppliers are denoted by $p_t\in\reals_+^{n_s}$, the fixed prices at which goods
can be sold to consumers are denoted by $r\in\reals_+^{n_c}$, and the consumer
demand is $d_t\in\reals_+^{n_c}$. The prices and demand are random and
independent between time points, but are known at time $t$ for planning. The
inputs are $b_t\in\reals_+^{n_s}$, amounts bought from the suppliers,
$s_t\in\reals_+^{n_c}$, the amounts sold to the consumers, and
$z_t\in\reals_+^{m-n_s-n_c}$, the amounts transported across inter-node links.
The dynamics are given by
\[
h_{t+1} = h_t + \left(A^{\text{in}} - A^{\text{out}}\right)(b_t,s_t,z_t),
\]
where $A^{\text{in}}, A^{\text{out}} \in\reals^{n\times m}$; $A^{\text{in
(out)}}_{ij}$ is 1 if link $j$ enters (exits) node $i$ and 0 otherwise.

The dynamics may be expressed as a random linear dynamical system with augmented
state $x_t = (h_t, p_t, d_t)$, input $u_t = (b_t, s_t, z_t)$, and dynamics
matrices
\[
A_t = \left[ \begin{array}{ccc}
I & 0 & 0 \\
0 & 0 & 0 \\
0 & 0 & 0 \\
\end{array}\right],
\qquad
B_t = \left[\begin{array}{c}A^\mathrm{in} - A^\mathrm{out} \\ 
0 \\ 0 \end{array}
\right],
\qquad
c_t = \left[
    \begin{array}{c}
        0 \\
        p_{t+1} \\
        d_{t+1}
    \end{array}
\right],
\]
such that $x_t\in\reals^n$ with $n=\tilde n + n_s + n_c$ and $u_t\in\reals^m$.

The prices and demand $p_t$ and $d_t$ are included in the state since they are
known at time $t$ for planning. However, since they are random and independent
between time points, the value function need only be a function of $h_t$.
Moreover, we only require that the stage cost be jointly convex in $(h_t, u_t)$.

The goal is to maximize the revenue from selling goods to customers while
minimizing the material costs paid to the suppliers, transportation costs, and
holding costs of the goods at each node. Let $\tau\in\reals_+^m$ encode the
costs of transporting a unit of good across each link, and $\alpha\in\reals^n_+$
and $\beta\in\reals^n_+$ parametrize the linear and quadratic holding costs of
the goods at each node.

The stage cost is
\[
g(x_t,u_t) = -r^T s_t + p_t^T b_t + \tau^T z_t 
+ \alpha^T h_t + \beta^T h_t^2 + I(x_t,u_t),
\]
where $I(x_t,u_t)$ is the indicator function that encodes the following
constraints:
\begin{itemize}
\item The warehouses have maximum capacity $h_{\max}>0$: $0 \le h_{t+1} \le
h_{\max}$. 
\item The links have maximum capacity $u_{\max}>0$: $0\le u_t \le u_{\max}$.
\item The amounts shipped out should not exceed the current capacities:
$A^\textrm{out}u_t \le h_t$.
\item The amounts sold to consumers cannot exceed the current demand: $s_t\le
d_t$.
\end{itemize}

For this example, we find a quadratic lower bound $V^{\textrm{lb}}$ on the value
function by relaxing the constraints, adding the quadratic penalty $u_t^T u_t -
(1/2)u_{\max}\ones^\top u_t$ to the stage cost, and solving the resulting LQR
problem. The lower bound is valid, since the added penalty is a pointwise lower
bound on the indicator of the input constraints, which is zero for $0 \le u_t
\le u_{\max}$, and infinity otherwise.

\PAR{Numerical example}
We consider a network with $\tilde n=4$ warehouses, $n_s=2$ suppliers, $n_c=2$
consumers, and $m=8$ links. The network is illustrated in Figure
\ref{f-sc-network}. The supplier prices $p_t$ and customer demands $d_t$ are
log-normally distributed, such that $\log p_t \sim
\mathcal{N}(\mu_p,\Sigma_p)$ and $\log d_t \sim \mathcal{N}(\mu_d,\Sigma_d)$,
with
\[
\mu_p = (0.0, 0.1), \quad \Sigma_p = 0.4I, \quad
\mu_d = (0.0, 0.4), \quad \Sigma_d = 0.4I.
\]
The holding cost parameters are $\alpha = \beta = (0.01)\ones$, the
transportation cost is $\tau = (0.05)\ones$, and the consumer prices are
$r=(1.3)\ones$. The maximum warehouse capacities are $h_{\max}=(3)\ones$, and
the maximum link capacities are $u_{\max}=(2)\ones$.

\PAR{Results}
We carried out VGI for $20$ iterations, starting from the quadratic lower bound
$V^{\mathrm{lb}}$.
In each iteration, the fitting step was performed using
$N=50$ fitting points, obtained by simulating the current policy. The damping
coefficient was fixed to $\rho_k=0.5$. When solving the fitting problem, we add
an $\ell_2$ (or ridge) regularization, with coefficient $\lambda = 10^{-4}$.

\begin{figure}
\centering
\includegraphics[scale=0.6]{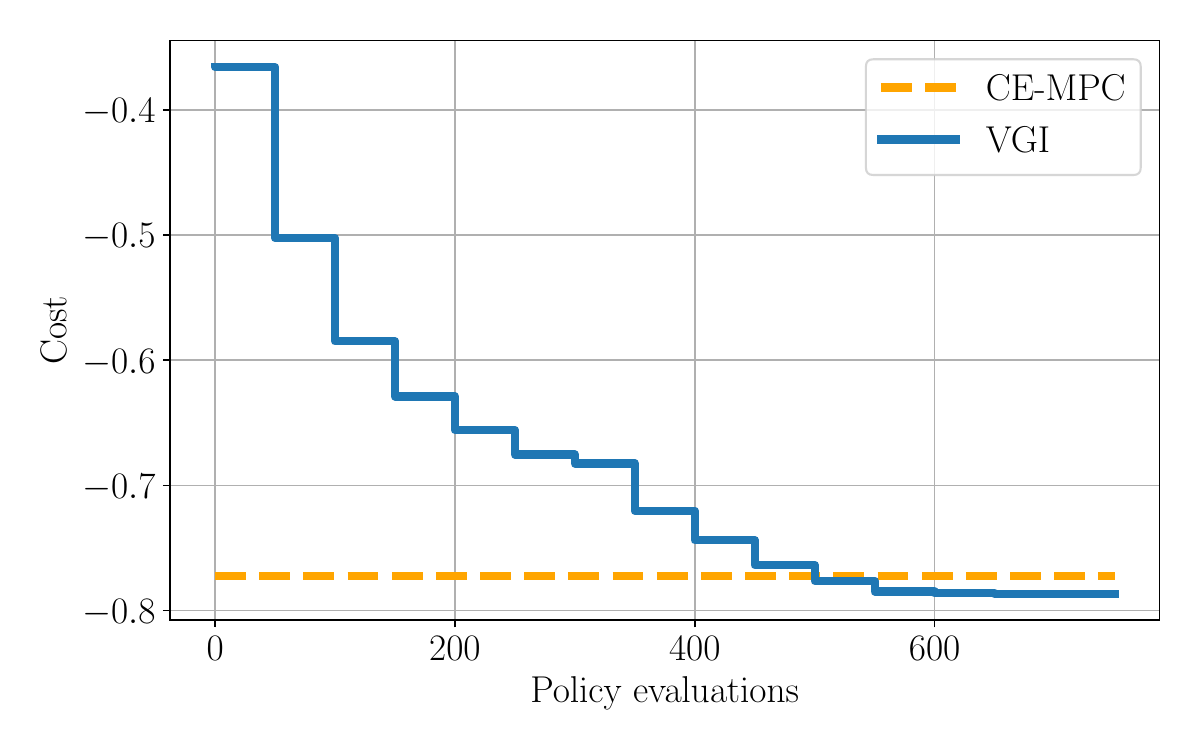}
\caption{VGI for the supply chain problem.}
\label{f-supply-chain-vgi}
\end{figure}

Figure \ref{f-supply-chain-vgi} shows the average cost versus the number of
policy evaluations used, along with the average cost of the CE-MPC policy with a
horizon of $H=30$. Our method converges to roughly the same cost as the CE-MPC
policy.

\begin{figure}
\centering
\includegraphics[scale=0.8]{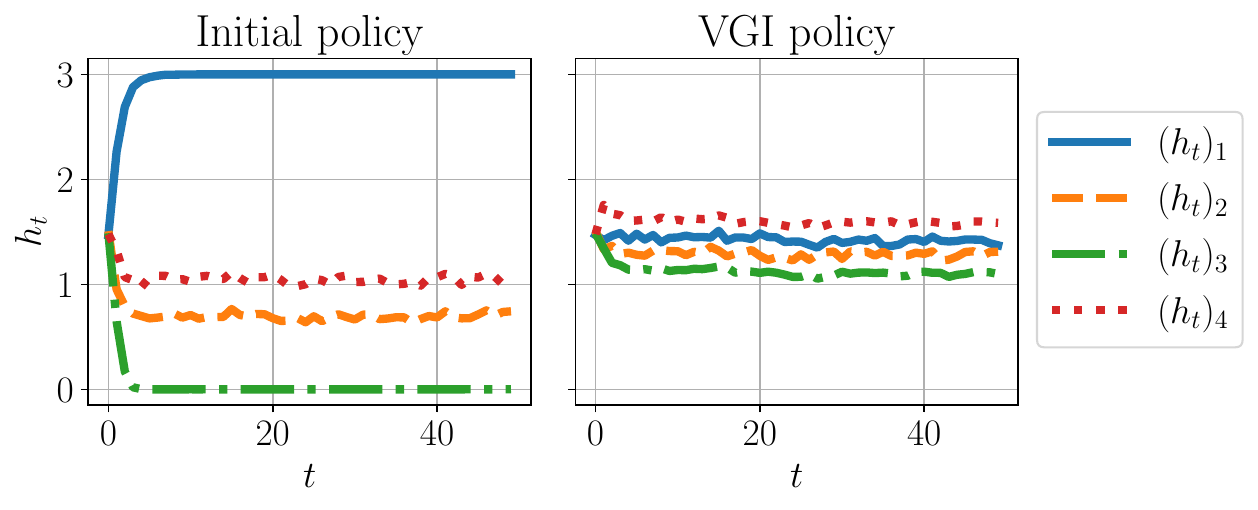}
\caption{Supply chain storage $h_t$ for each warehouse over time. Left: initial
ADP policy using $V^{\textrm{lb}}$. Right: final policy after VGI.}
\label{f-supply-chain-trajectories}
\end{figure}

Figure \ref{f-supply-chain-trajectories} shows the storage $h_t$ for each of the
four warehouses over time, for the initial policy using $V^{\textrm{lb}}$ and
the final policy after VGI. The plots show average trajectories over 500
simulations, each initialized with a state in $[0,h_{\max}]^4$, chosen uniformly
at random.

On average, the VGI policy is able to keep the storage levels close to half
capacity for all warehouses. On the other hand, the initial policy tends to put
too much stock in the first warehouse with storage $(h_t)_1$, which can, on
average, buy goods at a lower price from the suppliers. Similarly, the policy
tends to under-utilize the third warehouse with storage $(h_t)_3$, which
experiences lower consumer demand than the fourth warehouse with storage
$(h_t)_4$.

\section{Comparison with other methods}\label{s-comparisons}

In this section, we evaluate VGI against two related ADP methods for finding a
quadratic approximate value function: the standard FVI described in
\S\ref{s-fvi} and a COCP gradient method. They are iterative methods that follow
the same pattern as VGI: at each iteration, we simulate the system for $N$
steps, and then use the resulting data to update the parameters of the quadratic
approximate value function.

\PAR{COCP gradient method}
We compare against a gradient based method that updates the parameters $\theta$
of the ADP policy \eqref{e-adp-policy} using the derivatives of the cost along
simulated trajectories, with respect to $\theta$ \cite{AgrawalBBS2020}. At
iteration $k$, the policy $\phi^k$ with parameters $\theta^k$ is used to
simulate the system for $N$ steps. The resulting data is used to compute an
estimate of the average cost, given by
\[
\hat J(\theta^k) =
\frac{1}{N} \sum_{j=0}^{N-1} g(x_j, \phi^k(x_j)).
\]
We then compute $\nabla \hat J(\theta^k)$ using the chain rule, and then update
the parameters. This approach is known as backpropagation through time
\cite{Werbos1990}. In our experiments, we use the projected stochastic
(sub)gradient rule $\theta^{k+1} = \Pi_{\Theta}(\theta^k - \alpha^k \nabla\hat
J(\theta^k))$, where $\Pi_{\Theta}$ is the projection onto $\Theta$, and
$\alpha^k>0$ is a step size. 

This approach requires derivatives of the policy with respect to its parameters.
Those derivatives may be found by applying the implicit function theorem to the
optimality conditions of the convex optimization problem associated with the
policy \cite{AgrawalABBDK2019, AgrawalBBS2020}. Examples of the COCP gradient
method used to find quadratic approximate value functions may be found in
\cite{AgrawalBBS2020}. In our experiments, we used \texttt{cvxpylayers} to
compute the necessary derivatives \cite{AgrawalABBDK2019}.

\subsection{Results}

In general, FVI and COCP gradient methods required more tuning of
hyperparameters than VGI to work well. As shown in table \ref{t-comparison}, VGI
achieves the best (or close to the best) performance in all three problems, all
using far fewer policy evaluations than the FVI and COCP gradient methods. The
costs were evaluated in each case by simulating the policy for ten thousand
steps.

VGI used the same hyperparameters as in \S\ref{s-examples}, \ie, $\rho_k=0.5$
and $N=50$. The method was run for 40 iterations for the box-constrained LQR
problem, 20 iterations for the commitments example, and 15 iterations for the
supply chain problem.

We now discuss the hyperparameters chosen for FVI and the COCP gradient method.
All methods were initialized using the same initial quadratic approximate value
function. For the box-constrained LQR problem we used $V^1(x)=x^TQx$, and for
the other two problems we used $V^1(x)=V^{\textrm{lb}}$, the quadratic lower
bound on $V$ available for each problem.

\begin{table}
\centering
\caption{Comparison of cost and number of policy evaluations used (in thousands)}
\label{t-comparison}
\begin{tabular}{|c|cc|cc|cc|}
\hline
\multicolumn{1}{|c|}{} & \multicolumn{2}{c|}{\textbf{Box LQR}} &
\multicolumn{2}{c|}{\textbf{Commitments}} &
\multicolumn{2}{c|}{\textbf{Supply chain}} \\[5pt]
\multicolumn{1}{|c|}{\textbf{Method}} & 
\textbf{cost} & \textbf{evals. ($\times 10^3$)} & 
\textbf{cost} & \textbf{evals. ($\times 10^3$)} & 
\textbf{cost} & \textbf{evals. ($\times 10^3$)} \\
\hline
VGI           & 32.3 & 2 & 9.1 & 1 & -0.79 & 0.75 \\
FVI           & 32.2 & 20 & 9.1 & 4 & -0.77 & 16 \\
COCP gradient & 33.2 & 24 & 9.4 & 20 & -0.77 & 70 \\
MPC           & 33.3 & - & 11.9 & - & -0.77 & - \\
\hline
\end{tabular}
    
\end{table}
    
\PAR{Box constrained LQR} FVI was run using $N=400$ policy evaluations, for a
total of 50 iterations. The damping parameter was $\rho_k = 0.5$, and the
symmetry constraint $p=0$ was incorporated into the fitting problem.

The COCP gradient method was run using $N=300$ policy evaluations, for a total
of 80 iterations. The $300$ sample points were generated by simulating $K=3$
trajectories each of length $T=100$, using the procedure described in
\S\ref{s-parallel-sim}. We used a step size of $\alpha^k = 0.01$. The method was
initialized with $P=I$, and the symmetry constraint $p=0$ was incorporated into
the fitting problem. VGI took 6 seconds to complete, FVI took 29 seconds, and
the COCP gradient method took 4 minutes and 10 seconds.

\PAR{Commitments planning}
FVI was run using $N=200$ policy evaluations, for a total of $20$ iterations.
The sample points were generated by simulating $K=2$ trajectories each of length
$T=100$. The damping parameter was $\rho_k = 0.5$.

The COCP gradient method was run using $N=200$ policy evaluations, for a total
of 100 iterations. The sample points were generated by simulating $K=2$
trajectories each of length $T=100$. We used a step size of $\alpha^k =
10^{-4}$. VGI took 5 seconds to complete, FVI took 7 seconds, and the COCP
gradient method took 5 minutes.

\PAR{Supply chain}
FVI was run using $N=800$ policy evaluations, for a total of $20$ iterations.
The sample points were generated by simulating $K=2$ trajectories each of length
$T=400$. The damping parameter was $\rho_k = 0.75$.
An $\ell_2$ regularization with coefficient $\lambda=10^{-4}$ was used in the
fitting problem.

The COCP gradient method was run using $N=1000$ policy evaluations, for a total
of $70$ iterations. The sample points were generated by simulating $K=10$
trajectories each of length $T=100$. We used a step size of $\alpha^k = 0.01$.
An $\ell_2$ regularization with coefficient $\lambda=10^{-4}$ was added to the
cost. VGI took 2 seconds to complete, FVI took 25 seconds, and the COCP gradient
method took 13 minutes.

\section{Conclusion}\label{s-conclusion}

In this work, we propose value-gradient iteration, a method for finding a
quadratic approximate value function for convex stochastic control. The method
is an approximation of value iteration, and we show how we may compute the
gradient of the Bellman operator image to fit the gradient of the approximate
value function in each iteration. By fitting the gradient of the approximate
value function instead of the approximate value function itself, we can find a
good policy using far less simulation data. Indeed, we find that the
computational effort of obtaining a good approximate value function is
comparable to that of evaluating the policy through simulation.

\subsection*{Acknowledgements}

Stephen Boyd was partially supported by ACCESS (AI Chip Center for Emerging
Smart Systems), sponsored by InnoHK funding, Hong Kong SAR, and by Office of
Naval Research grant N00014-22-1-2121. We would also like to thank Pieter Abeel
and Zico Kolter for helpful discussions and feedback.

\clearpage
\bibliography{vgi}
\clearpage

\appendix

\section{Expectation of quadratic functions}\label{a-expect-quad-fns}

Let $\hat V$ be a convex quadratic function of the form
\eqref{e-quad-value-function}. We now show that $\Expect\hat V(A_tx+B_tu+c_t)$
is a convex quadratic function in $u$, with coefficients that may be written in
terms of $P$, $p$, and $\pi$ and the first and second moments of
$(A_t,B_t,c_t)$.

Let $\bar A = \Expect A_t$, $\bar A = \Expect A_t$, and $\bar c = \Expect c_t$
denote the expected values. Let $(A_t)_i$, $(B_t)_i$, $\bar A_i$, and $\bar B_i$
denote the $i$th columns of $A_t$, $B_t$, $\bar A$, and $\bar B$ respectively.
Let $\Sigma^{AB}_{ij}$ denote the covariance matrix between the $i$th column of
$A_t$ and the $j$th column of $B_t$, and let $\Sigma^{A}_{ij}$ and
$\Sigma^{B}_{ij}$ be defined similarly. Finally, let $\Sigma^{Ac}_{i}$ and
$\Sigma^{Bc}_{i}$ denote the covariances between the $i$th columns of $A_t$ and
$B_t$ with $c_t$, respectively, and let $\Sigma^c$ denote the covariance of
$c_t$.

We have
\[
\Expect\hat V(A_tx+B_tu+c_t) = \frac 1 2 \Expect\left(
\BBM A_tx+B_tu+c_t \\ 1 \EBM^T
\BBM P & p \\ p^T & \pi \EBM
\BBM A_tx+B_tu+c_t \\ 1 \EBM\right).
\]
Expanding terms, we obtain
\[
\Expect\hat V(A_tx+B_tu+c_t) = \frac 1 2
\BBM u \\ 1 \EBM^T
\BBM M & m \\ m^T & \mu\EBM
\BBM u \\ 1 \EBM,
\]
where
\[
\begin{array}{ll}
M = \Expect B_t^TPB_t, \\
m = \Expect B_t^TPA_tx + \Expect B_t^TPc_t + \bar B^T p, \\
\mu = \pi + x^T\Expect(A_t^TPA_t)x + 2x^T\Expect(A_t^TPc_t) + 2x^T\bar A^Tp
+ 2p^T\bar c + \Expect c_t^TPc_t.
\end{array}
\]
Finally, we note that
\[
\Expect c_t^T P c_t = \bar c^T P \bar c + \Tr(P\Sigma^{c}),
\]
and for all indices $i$ and $j$,
\[
\begin{array}{ll}
(\Expect A_t^T PA_t)_{ij} = \bar A_i^T P\bar A_j^T + \Tr(P\Sigma^{A}_{ij}),
\\[4pt]
(\Expect B_t^T PB_t)_{ij} = \bar B_i^T P\bar B_j^T + \Tr(P\Sigma^{B}_{ij}),
\\[4pt]
(\Expect B_t^T PA_t)_{ij} = \bar B_i^T P\bar A_j^T + \Tr(P\Sigma^{AB}_{ij}),
\\[4pt]
(\Expect A_t^T Pc_t)_i = \bar A_i^TP\bar c + \Tr(P\Sigma^{Ac}_i),
\\[4pt]
(\Expect B_t^T Pc_t)_i = \bar B_i^TP\bar c + \Tr(P\Sigma^{Bc}_i).
\end{array}
\]

\section{Lower bounds on quadratic functions}\label{a-quad-lower-bound}

We say that $V_1\ge V_2$ if $V_1(x)\ge V_2(x)$ for all $x\in\reals^n$. We
consider the case of convex quadratic functions, where for $i=1,2$, $V_i$ is
given by
\[
V_i(x) = 
\frac 1 2
\BBM
x \cr 1
\EBM^T
\BBM
P_i & p_i \cr
p_i^T & \pi_i
\EBM
\BBM
x \cr 1
\EBM,
\]
where $P_i\succeq 0$. Then, $V_1\ge V_2$ holds if and only if the quadratic
function $V_{12}=V_1-V_2$ is positive semidefinite, \ie
\[
V_{12}(x) = 
\frac 1 2
\BBM
x \cr 1
\EBM^T
\BBM
P_1 - P_2 & p_1-p_2 \cr
p_1^T-p_2^T & \pi_1-\pi_2
\EBM
\BBM
x \cr 1
\EBM \ge 0,
\]
for all $x\in\reals^n$. The function $V_{12}$ has a minimum value if and only if
$P_1-P_2\succeq 0$ and $p_1-p_2$ is in the range of the matrix $P_1-P_2$ (see
\eg~\cite[\S A.5.5]{BoydV2004}). The range condition may be written as
\[
\left[I - (P_1-P_2)(P_1-P_2)^\dagger\right](p_1 - p_2) = 0,
\]
where $(P_1-P_2)^\dagger$ is the pseudo-inverse of $(P_1-P_2)$. In this case,
the minimum value is given by
\[
\min_x V_{12}(x) = \frac 1 2 \left(\pi_1-\pi_2 - (p_1-p_2)^T
(P_1-P_2)^\dagger(p_1-p_2)\right).
\]
Finally, by the generalized Schur complement, $\min_x V_{12}(x)$ exists and is
nonnegative if and only if
\[
\BBM
P_1 - P_2 & p_1 - p_2 \cr 
p_1^T - p_2^T & \pi_1-\pi_2
\EBM\succeq 0.
\]

\section{Lower bound from certainty equivalence}\label{a-ce-lb}

Solving the certainty equivalent problem involves finding a function
$V^{\text{ce}}$ that satisfies the Bellman equation
\[
V^{\text{ce}}(x) = \min_u \left(g(x,u) 
+ V^{\text{ce}}(\bar Ax + \bar Bu + \bar c)\right).
\]
By Jensen's inequality,
\[
V^{\text{ce}}(\bar Ax+\bar Bu+\bar c)\le\Expect V^{\text{ce}}(Ax+Bu+c).
\]
Therefore, we have
\[
V^{\text{ce}}(x) \le \min_u \left(g(x,u) 
+ \Expect V^{\text{ce}}(Ax + Bu + c)\right) = \mathcal T V^{\text{ce}}(x).
\]
By the monotonicity of the Bellman operator, we have
\[
V^{\text{ce}} \le \mathcal T V^{\text{ce}}
\le \lim_{k\to\infty} \mathcal T^k V^{\text{ce}} = V^\star.
\]
This implies that $V^{\text{ce}}$ is a lower bound on the true value function.

\end{document}